\documentclass[journal,letterpaper,10pt]{IEEEtran}

\usepackage[cmex10]{amsmath}
\usepackage{algorithm}
\usepackage{algpseudocode} 
\usepackage{graphicx}
\usepackage{amsmath,amssymb}
\usepackage{color}

\newcommand{\Ne}[1]{\ensuremath{N(#1)}}
\newcommand{\Nhminus}[1]{\ensuremath{H^{-}_{#1}(t)}}
\newcommand{\Nhminusw}[1]{\ensuremath{\hat{H}^{-}_{#1}(t)}}
\newcommand{\Nsminus}[1]{\ensuremath{S^{-}_{#1}(t)}}

\newcommand{\Merr}[1]{\ensuremath{M\left(#1\right)}}
\newcommand{\Terr}[1]{\ensuremath{T\left(#1\right)}}
\newcommand{\mR}{\mathbb{R}}
\newcommand{\expec}[1]{\textbf{E}\left[#1\right]}
\newcommand{\LSTV}{LS-TV}
\newcommand{\LSMV}{LS-MV}

\newcommand{\pos}[2]{\ensuremath{p_{#1#2}}}

\newcommand{\ssv}[1]{\ensuremath{\sigma_{#1}}}
\newcommand{\Deg}[1]{\ensuremath{\Delta_{#1}}}
\newcommand{\xt}{\ensuremath{\overline{x}}}
\newcommand{\Lff}{\ensuremath{L_{ff}}}
\newcommand{\Lfl}{\ensuremath{L_{fl}}}

\newcommand{\xf}{\ensuremath{x^\textit{f}(t)}}
\newcommand{\xl}{\ensuremath{x^l(t)}}

\newcommand{\tp}{\ensuremath{^{\mathsf{T}}}}
\newcommand{\leaderSet}{\ensuremath{U}}

\newcommand{\resist}[2]{\Omega_{#1#2}}
\newcommand{\dist}[2]{D_{#1#2}}
\newcommand{\fdist}[2]{d_{#1#2}}

\newcommand{\resistwt}[2]{\Omega_{#1#2}}
\newcommand{\Gwt}{\ensuremath{\hat{G}}}
\newcommand{\Ewt}{\ensuremath{\hat{E}}}

\newtheorem{theorem}{Theorem}
\newtheorem{definition}{Definition}

\author {Stacy Patterson
\thanks{S. Patterson is with the Department of Computer Science, Rensselaer Polytechnic Institute, Troy, NY 12180, USA, 
 {\tt\small sep@cs.rpi.edu}}%
 }

\title{\LARGE \bf In-Network Leader Selection for Acyclic Graphs}

\begin{document}
\maketitle
\thispagestyle{empty}
\pagestyle{empty}

%\title{\LARGE \bf Optimal, In-Network Leader Selection for Acyclic Graphs}

\begin{abstract}
%\bold math
We study the problem of leader selection in leader-follower multi-agent systems that are subject to stochastic disturbances.
This problem arises in applications such as vehicle formation control, distributed clock synchronization, and distributed localization in sensor networks.
We pose a new leader selection problem called the in-network leader selection problem.  
Initially, an arbitrary node is selected to be a leader, and in all consequent steps the network must have exactly one leader.
The agents must collaborate to find the leader that minimizes the variance of the deviation from the desired trajectory, and they must do so within the network using only communication between neighbors.  
To develop a solution for this problem, we first show a connection between the leader selection problem and a class of discrete facility location problems.
We then leverage a previously proposed self-stabilizing facility location algorithm to develop a self-stabilizing in-network leader selection algorithm for acyclic graphs.  
\end{abstract}

\vspace{-.2cm}

\section{Introduction}
%One paragraph on networked control systems.

% One paragraph on leader-follower systems.
We consider a class of multi-agent systems called \emph{leader-follower systems}.
One or more agents in the network are \emph{leaders};  the states of these agents serve as the reference
states for the system.  The remaining agents are \emph{followers} that update their states based on relative information
about their own states and the states of their neighbors.  Hence, a system owner may control the entire network by controlling the leader agents.
These dynamics arise  in applications such as vehicle formation 
control~\cite{RBM05}, distributed clock synchronization~\cite{EKPS04}, and distributed localization in sensor networks~\cite{BH09}, for example.

It has been shown that the performance of leader-follower systems, where the followers are subject to stochastic disturbances,
 depends on the location of the leaders in the network~\cite{BH06,PB10}. This relationship naturally raises the question of how to choose leaders, from among all agents in the network, that give the best performance.  The \emph{$k$-leader selection problem} posed in~\cite{PB10},  is to select $k$ leaders  that minimizes the \emph{network coherence}, an $H_2$ norm of the leader-follower system that 
quantifies the variance of the nodes states from the target states.    The optimal leader set can be found by an exhaustive search over all subsets of size $k$, but this solution is not tractable in large networks.  

Several recent works have proposed efficient approximation  algorithms for the $k$-leader selection problem. 
In these works, the full network topology is known, and the leader set is selected using an offline algorithm.
In~\cite{LFJ14}, 
the authors use a convex relaxation of the combinatorial $k$-leader selection problem, and they present an efficient interior point method this relaxed problem.    In~\cite{CBP14}, the authors show that the mean-square deviation from the desired state is proportional to super-modular set function~\cite{N78}.
As such, a greedy, polynomial time solution can be used to find a leader set for which the mean-square error is within a provable bound of optimal.
Other works have explored the optimal leader selection in leader-follower systems without stochastic disturbances~\cite{CABP14} and in systems where both
the leaders and followers are subject to stochastic disturbances and the leaders also have access to relative state information~\cite{LFJ14}. 
Finally, recent work has shown that the optimal single leader and optimal pair of leaders in a network are those nodes with maximal information centrality and joint centrality, respectively~\cite{FL13}.

In this work, we investigate a variation on the $k$-leader selection problem for $k=1$ that we call the \emph{in-network leader selection problem}.
We consider two performance measures, the total variance of the deviation from the desired trajectory and the maximum variance of this 
deviation, over all agents.
Initially, a single agent is selected as the leader, and the network must have a single leader at all times.  The agents must collaborate to find the leader that minimizes the specified performance measure (total or maximal variance),  and they must do so within the network using only communication between neighbors.  This problem may arise in a multi-agent system that has limited bandwidth to the system owner.
The system owner controls the network through a single communication channel with the leader. 
The system can determine the optimal leader within the network
and inform the owner when leadership is transferred from one agent to another.  Ideally, should the network topology change, the leadership should be transferred to the new optimal leader.

We first show a connection between optimal leader selection and two discrete facility location problems, the
$p$-median problem and the $p$-center problem~\cite{HD02}.  We then leverage previously proposed algorithms for these facility location problems to develop a solution for the in-network leader selection problem for acyclic graphs.  Our approach is self-stabilizing, meaning that if, after the algorithm has found an optimal leader,   the network topology changes, the algorithm will then find an optimal leader for the new topology.

The remainder of this work is organized as follows.  In Section~\ref{model.sec},
we present the system model and background on the leader selection problem. In Section~\ref{problem.sec}, we formalize the in-network
leader selection problem.  Section~\ref{algorithm.sec} gives our algorithm and its analysis, including the relationship to the $p$-median 
and $p$-center facility location problems.  In Section~\ref{extend.sec}, we show how our algorithm can be extended to leader-follower systems
in weighted graphs.  Finally, we conclude in Section~\ref{conclusion.sec}.

%%%%%%%%%%%%%%%%%%%%%%%%%%%%%%%%%%%%%%%%%%%%%%%%%%%%%%%%%%%%%%%%%%%%%%%%%
%%%%%%%%%%%%%%%%%%%%%%%%%%%%%%%%%%%%%%%%%%%%%%%%%%%%%%%%%%%%%%%%%%%%%%%%%
%%%%%%%%%%%%%%%%%%%%%%%%%%%%%%%%%%%%%%%%%%%%%%%%%%%%%%%%%%%%%%%%%%%%%%%%%
%%%%%%%%%%%%%%%%%%%%%%%%%%%%%%%%%%%%%%%%%%%%%%%%%%%%%%%%%%%%%%%%%%%%%%%%%
%%%%%%%%%%%%%%%%%%%%%%%%%%%%%%%%%%%%%%%%%%%%%%%%%%%%%%%%%%%%%%%%%%%%%%%%%

\vspace{-.2cm}

\section{Background and System Model}
\label{model.sec}

We consider a system of $n$ interconnected agents or nodes.  The communication structure is modeled by a connected, undirected graph
$G = (V,E)$, where $V$ is the set of of agents, with $|V| = n$ and $E$ is the set of communication links, with $|E| = m$. An edge $(i,j)$ exists in $E$ if and only if node $i$
and node $j$ can send information between them.  We denote the neighbor set of a node $i$ by $\Ne{i}$.

Every agent has a scalar-valued state $x_i(t)$, and the state of the system is given by the vector $x(t) \in \mR^n$.
A subset of the agents $\leaderSet \subseteq V$ are \emph{leaders}.  We assume that the states of these agents serve as a reference for the network;
these states remain fixed at an identical, constant trajectory $\xt \in \mR$, i.e.,
\begin{equation}
x_i(t) = \xt,~\text{for all}~ {t \geq 0},~\text{for all}~i \in \leaderSet. \label{leader.eq}
\end{equation}
%Without loss of generality, we assume that $\xt = 0$ (see \cite{PB10}).

The remaining agents, those in $V~\backslash~\leaderSet$, are \emph{followers}.  
Each follower $i$ updates its state based on its own state and those of its neighbors using a simple noisy consensus algorithm,
\begin{equation}
\dot{x}_i(t) = \sum_{j \in \Ne{i}} x_i(t) - x_j(t)  + w_i(t), \label{follower.eq}
\end{equation}
where $w_i(t)$ is a white stochastic disturbance with zero-mean and unit variance.

Without loss of generality, we let the nodes be ordered so that $x(t) = [\xf\tp~\xl\tp]\tp$, where $\xf$ denotes is the $|V \setminus \leaderSet|$-vector
containing the states of the follower agents and $\xl$ is the $|\leaderSet|$-vector containing the states of the leader agents.
Let $L$ be the Laplacian matrix of the graph induced by the leader-follower dynamics.  For now, we restrict our study to consensus dynamics over an unweighted graph, and thus we use the unweighted Laplacian matrix.    In Section~\ref{extend.sec}, we show how our problem formulation 
and solution can be extended to dynamics that depend on a weighed Laplacian matrix.
Each component of $L$ is defined as
\[
L_{ij} = \left(\begin{array}{ll}
-1&~\text{if}~(i,j) \in E~\text{and $i$ is a follower} \\
\text{degree($i$)}&~\text{$i=j$ and $i$ is a follower} \\
0&~\text{otherwise.}
\end{array} \right.
\]
$L$ can be decomposed according to the leader/follower designations as
\[
L = \left[ \begin{array}{cc}
\Lff & \Lfl \\
0 & 0 
\end{array} \right],
\]
where $\Lff$ defines the interactions between followers and $\Lfl$ defines the impact of the leaders on the followers.
We note that (since $G$ is connected), $\Lff$ is positive definite~\cite{PB10}.

With this decomposition, the evolution of the states of the follower nodes can be expressed as,
\begin{equation}
\dot{x}^f(t) = - \Lff x(t) + w^f(t),
\end{equation}
where $w^f(t)$ is the $|V \setminus \leaderSet|$-vector of stochastic disturbances.

In the absence of the noise terms $w_i(t)$, the agents would converge to the desired state $\xt$.  With these noise terms, the
agents' states do not converge, however, the steady-state variance of their deviations from $\xt$ are bounded~(see e.g., \cite{PB10}).
Formally, we define the steady-state variance of agent $i$ as
\[
\ssv{i} = \lim_{t \rightarrow \infty} \expec{(x_i(t) - \xt)^2}.
\]
It was shown in ~\cite{BH06} that this variance is related to the $(i,i)^{th}$ entry of $\Lff^{-1}$ as $\ssv{i} = \frac{1}{2} (\Lff^{-1})_{ii}$.
%\begin{theorem}[\cite{BH05}, Section IV-B]
%For a given leader set $S$, let $\Lff$ denote the corresponding follower sub-matrix of the Laplacian.  The steady-state variance of the deviation from $\xt$ at an agent $i$ is given by 
%\[
%\ssv{i} = \frac{1}{2} (\Lff^{-1})_{ii}.
%\]
%\end{theorem}

It has also been shown that, for $|\leaderSet| = 1$, the steady-state variance of a node $i$ can also be expressed in terms of the 
\emph{resistance distance} from $i$ to the leader, where the resistance distance is defined as follows.
Let the graph represent an electrical network where each edge is a unit resistor.
The resistance distance $\resist{i}{j}$ between two nodes $i$ and $j$ is the potential difference between
them when a one ampere current source is connected from node $j$ to node $i$.
\begin{theorem}[See \cite{GBS08}]
For a single leader $s$, the steady-state variance of the deviation from $\xt$ at an agent $i$  is related to the resistance distance between $s$ to $i$ as
 $\ssv{i} = \frac{1}{2} \resist{i}{s}$.
\end{theorem}

 In a general graph, resistance distance is a distance function.  For a simple, connected graph only, the resistance distance between nodes $i$ and $j$ is equivalent to the  
 conventional graph distance $\dist{i}{j}$ where $\dist{i}{j}$ is the length of the shortest path between $i$ and $j$~\cite{KR93}.

\section{Problem Formulation}
\label{problem.sec}

The steady-state variance of each agent depends on the choice of the leader set $\leaderSet$.  
The leader selection problem involves identifying a set $\leaderSet$
that minimizes a function of these variances. We next define  the functions that we use to measure the performance of a leader set, followed 
by a formal definition of the leader selection problems we address in this work.

\subsection{Performance Measures}
We quantify the relationship between the leader set and the steady-state variance using two performance measures.
The first is the \emph{total steady-state variance}, which is,
\[
\Terr{\leaderSet} := \sum_{i \in V \setminus U} \ssv{i}.
\]
This error measure, related to the coherence of the network~\cite{BJMP12}, has been studied in previous works on leader selection~\cite{PB10,LFJ14,CABP14}.
We also consider the \emph{maximum steady-state variance}
over all agents, 
\[
\Merr{\leaderSet} := \max_{i \in V \setminus U} \ssv{i}.
\]
As far as we are aware, this error measure has not been considered in previous works.

\begin{figure}
\centering
\includegraphics[scale=.6]{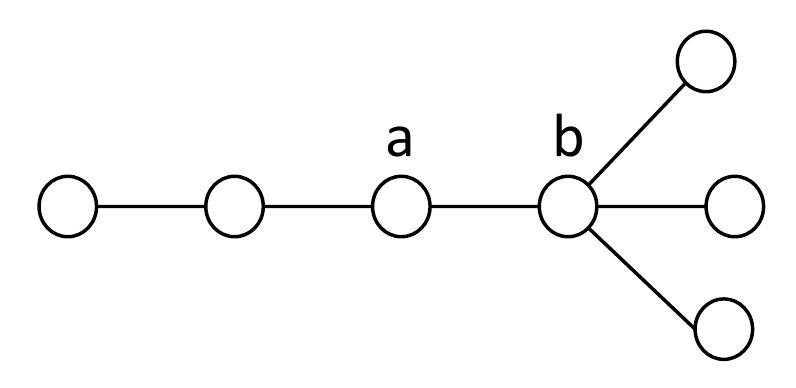}
\caption{Example graph with optimal leaders.  For \LSTV, the optimal leader is node $b$,
for which $\Terr{\{b\}}  =  4.5$ and $\Merr{\{b\}} = 1.5$.  For \LSMV, the optimal leader is node $a$, for which $\Terr{\{a\}} = 5$ and $\Merr{\{a\}} = 1$.}
\label{mixed_graph.fig}
\end{figure}

\subsection{The Leader Selection Problem}

The goal of the leader selection problem is to identify a leader set $\leaderSet$ such that the steady-state variance of the agents in minimized.
The \emph{total variance $k$-leader selection problem} ($k$-\LSTV) is
\begin{equation} \label{Tleaderprob.eq}
\begin{array}{cc}
\text{minimize} &~\Terr{\leaderSet}\\
\text{subject to}&~|\leaderSet| = k.
\end{array}
\end{equation}
The \emph{maximal variance $k$-leader selection problem} ($k$-\LSMV) is 
\begin{equation} \label{Mleaderprob.eq}
\begin{array}{cc}
\text{minimize} &~\Merr{\leaderSet}\\
\text{subject to}&~|\leaderSet| = k.
\end{array}
\end{equation}
For $k=1$, we omit $k$ from the naming convection, denoting these problems by \LSTV\ and \LSMV.
We note that the optimal leader set may be different depending on which performance measure is used.
An example of this is given in Figure~\ref{mixed_graph.fig} for $k=1$.

A naive solution to the this problem is to compute $\Terr{\leaderSet}$ ($\Merr{\leaderSet}$) for all subsets $\leaderSet \subseteq V$ with $|\leaderSet| = k$ and to choose the leader set
for which this function is minimized.  However, this approach has combinatorial complexity. Several works have proposed approximation algorithms for  \LSTV,
that run in polynomial time.  Notably, in~\cite{CBP14}, it was shown that $\Terr{\leaderSet}$ is proportional to a super-modular function, which implies that a simple greedy
(polynomial time) leader selection scheme yields a leader set whose error is within a bounded approximation of optimal.  We note that the error function 
$\Merr{\leaderSet}$ is not super-modular. A  proof this is given in Appendix~\ref{notsubmod.app}.  Therefore, the results for the approximation algorithm for \LSTV\ do
not necessarily extend to \LSMV.

\vspace{-.2cm}

\subsection{The In-Network Leader Selection Problem}

We now propose a variation on the $k$-leader selection problem where agents collaborate to determine the optimal leader set. For this initial investigation of the problem, we restrict our focus to the case where $k=1$ and where the network is a connected,
acyclic graph.

Initially, a single agent is selected as leader arbitrarily.  The state of this leader is the reference state $\xt$ for the network.
In each round, every agent communicates with its neighbors.
During this communication, the leader agent may choose to transfer the leadership role from itself to one of its neighbors.  
Since the system is synchronous, the old and new leader can schedule the leadership handoff so that it occurs instantaneously.
The  new leader adopts the reference state $\xt$ and then follows the dynamics in (\ref{leader.eq}).
After an agent transfers its leadership, it behaves as a follower according to the dynamics in (\ref{follower.eq}). 
The dynamics of the remaining nodes remain unchanged after the handoff.
The goal is for the leadership role to eventually be transferred to and remain at the leader $u$ for which $F(\{u\})$ is minimized, where
$F(\cdot) = \Merr{\cdot}$ or $F(\cdot) = \Terr{\cdot}$.

Our aim is to develop algorithms for in-network leader selection that are \emph{self-stabilizing}, which is formally defined as follows.
\begin{definition} \label{stable.def}
A leader selection algorithm is \emph{self-stabilizing} if, starting with any initial leader, the leadership role is transferred to the optimal leader in a finite number of steps (\emph{convergence}) \emph{and}
it remains there for all subsequent steps (\emph{closure}). \end{definition}
A self-stabilizing algorithm is robust to changes in the network.
For example, suppose the optimal leader is selected, and then the network topology changes.
Provided that after this change, the network topology remains stable for ``long enough'', the  algorithm will then find the optimal leader for the new topology. 

The standard definition of a self-stabilizing algorithm~\cite{D74} requires that the algorithm converge to the desired solution from any initial state.
Definition~\ref{stable.def} differs from this definition slightly in that we require that a single leader is selected in the initial state.
This assumption implies that a self-stabilizing leader selection algorithm need not be robust to agent failures since the failure of the leader agent
violates the assumption.  Conceivably, under certain assumptions about the frequency of failures, a distributed leader election algorithm~\cite{L96} could be used to replace a failed leader agent.  The integration of leader election and optimal in-network leader selection is a subject for future work.

%%%%%%%%%%%%%%%%%%%%%%%%%%%%%%%%%%%%%%%%%%%%%%%%%%%%%%%%%%%%%%%%%%%%%%%%%
%%%%%%%%%%%%%%%%%%%%%%%%%%%%%%%%%%%%%%%%%%%%%%%%%%%%%%%%%%%%%%%%%%%%%%%%%
%%%%%%%%%%%%%%%%%%%%%%%%%%%%%%%%%%%%%%%%%%%%%%%%%%%%%%%%%%%%%%%%%%%%%%%%%
%%%%%%%%%%%%%%%%%%%%%%%%%%%%%%%%%%%%%%%%%%%%%%%%%%%%%%%%%%%%%%%%%%%%%%%%%
%%%%%%%%%%%%%%%%%%%%%%%%%%%%%%%%%%%%%%%%%%%%%%%%%%%%%%%%%%%%%%%%%%%%%%%%%

\vspace{-.2cm}

\section{In-Network Leader Selection Algorithm}
 \label{algorithm.sec}

In this section, we present our algorithms for in-network leader selection.
To develop our algorithms, we leverage connections between the leader selection problem and two discrete facility location problems,
the $p$-median problem and the $p$-center problems.  We next present a summary of these problems, followed by a description of our algorithms
and their analysis.
\vspace{-.2cm}

\subsection{The $p$-Median and $p$-Center Problems}

The $p$-median and $p$-center problems belong to a larger class of discrete facility location problems~\cite{HD02}.
In this class of problems, there is a discrete set of demand nodes, a discrete set of candidate facility locations,
and specified distance $\fdist{i}{j}$ between each demand node $i$ and candidate location $j$.   Each demand node will
be assigned to the closest facility.
The objective is to select a set of $p$ facility locations that minimizes some function of the distances between the demand
nodes and their assigned facilities.  

In both the $p$-median and $p$-center problems, the candidate facility locations coincide with the locations of the demand nodes.
We denote this set of nodes by $V$.
For the $p$-median problem, $p \leq |V|$ facility locations are selected so as to minimize the sum of the distances between each demand node and its assigned facility,
\begin{equation} \label{pmedianprob.eq}
\begin{array}{ll}
\text{minimize} &\displaystyle ~\sum_{i \in V} \min_{u \in \leaderSet} \fdist{i}{u} \\
\text{subject to} & \leaderSet \subseteq V,~ |\leaderSet| = p.
\end{array}
\end{equation}
In the $p$-center problem, $p$ facility locations are selected to minimize the maximum distance between any demand node and its assigned facility,
\begin{equation} \label{pcenterprob.eq}
\begin{array}{ll}
\text{minimize} &\displaystyle ~\max_{i \in V} \min_{u \in \leaderSet} \fdist{i}{u} \\
\text{subject to} &~ \leaderSet \subseteq V,~|\leaderSet| = p.
\end{array}
\end{equation}
In general, the $p$-median problem is NP-Hard, and the $p$-center problem is NP-complete.  In simple graphs, however, these problems can be solved in polynomial-time~\cite{KH79}.

If we consider the agents of the network as the set of demand nodes/candidate facility locations and let the $\fdist{i}{j} =  \resist{i}{j}$ for all nodes $i,j$,
then for $p = 1$, a solution to the $p$-median problem gives a solution to \LSTV.  Similarly, a solution to the $p$-center problem gives a solution to 
\LSMV.

\begin{algorithm}
\caption{Algorithm for in-network leader selection for \LSTV.} \label{leader.alg}
\begin{algorithmic}
\State {Send $s_i(t)$ to all  agents $j$ in $\Ne{i}$}
\State {Receive $s_j(t)$ from all $j \in \Ne{i}$}
\If {$|\Ne{i}| = 1$}
	\State{$s_i(t+1) \gets 1$}
\Else
	\State {$s_i(t+1)  \gets  1+ \sum \left(\Nsminus{i}\right)$}
\EndIf
\If {(\textit{leader} = TRUE) and \\~~~~~~~~~~~($s_j(t) > s_i(t+1)$ for some $j \in \Ne{i}$) }
	\State{$k =  \arg \max_{j \in \Ne{i}} s_j(t)$}
	\State{$\textit{leader} \gets$  FALSE}
	\State{transfer leadership to agent $k$}
\ElsIf{receive leadership transfer from neighbor}
	\State{$\textit{leader} \gets$ TRUE}
\EndIf
\end{algorithmic}
\end{algorithm}

\subsection{Self-Stabilizing Leader Selection Algorithm}

We now describe our in-network leader selection algorithm for acyclic graphs.  Recall that in an unweighted, acyclic, connected graph, the resistance distance $\resist{i}{j}$  between nodes $i$ and $j$  is equal to the graph distance $\dist{i}{j}$\footnote{The resistance distance between neighboring nodes $i$ and $j$ is $\resist{i}{j} = 1$.}.
When the distance $\fdist{i}{j}$ in (\ref{pmedianprob.eq}) and (\ref{pcenterprob.eq}) is given by the graph distance, a solution to the $1$-median problem is called a \emph{median} of the graph and a solution to the $1$-center problem is called a \emph{center} of the graph. Thus, a median of the graph is an optimal leader for \LSTV, and a center of the graph is an optimal leader for \LSMV.
A graph may have more than one median or center;
should the graph have multiple medians or centers, any one can be selected as the optimal leader.

Our approach for in-network leader selection is based on a simple self-stabilizing approach for finding the median and center of graph that was proposed in~\cite{BGKP99}, which we summarize below. For ease of presentation, we adopt a synchronized communication model.  Communication takes place in rounds, and in each round, an agent exchanges information with all of its neighbors. We note median and center-finding algorithms have been proven correct under much less restrictive communication models.

In the median finding algorithm, each agent has a variable $s_i(t)$.These variables are called the \emph{$s$-values} of the agents. 
Since the algorithm is self-stabilizing, there is no need to specify an initial value for $s_i(0)$; the algorithm
will converge to the correct solution from any initial value. In each round, the agent sends its $s_i(t)$ to all of its neighbors.  
Let $\Nsminus{i}$ be the set of values $s_i(t)$ received from $j \in \Ne{i}$ in round $t$, with one maximum $s_j(t)$ removed.
The agent then updates $s_i(t)$ as follows,
\begin{equation} \label{medianalg.eq}
s_i(t+1) = \left\{ \begin{array}{ll}
1 &~\text{if}~|\Ne{i}| = 1 \\
1 + \sum \left(\Nsminus{i}\right)&~\text{otherwise}.
\end{array}\right.
\end{equation}
Here $\sum \left(\Nsminus{i}\right)$ denotes the sum over the elements in the set $\Nsminus{i}$.

The center finding algorithm operates in a  similar manner.  Each agent has a variable $h_i(t)$. These variables are called the \emph{$h$-values} of the agents.  
In each round, all neighbors exchange their $h$-values.
Let $\Nhminus{i}$ be the set of values $h_i(t)$ received from $j \in \Ne{i}$, with one maximum $h_j(t)$ removed.
The agent updates $h_i(t)$ as follows,
\begin{equation} \label{centeralg.eq}
h_i(t+1) = \left\{ \begin{array}{ll}
0 &~\text{if}~|\Ne{i}| = 1 \\
1 + \max \left(\Nhminus{i}\right)&~\text{otherwise}.
\end{array}\right.
\end{equation}
Here $\max \left(\Nhminus{i}\right)$ denotes the maximal value in the set $\Nhminus{i}$.

The following theorem gives a formal statement of the convergence and closure guarantees of these algorithms 
\begin{theorem}[See \cite{BGKP99}] \label{selfstable.thm}
There exists a finite time $T$ such that $s_i(t+1) = s_i(t)$ 
($h_i(t+1) = h_i(t)$) for all $i \in V$, for all $t \geq T$.  For all $t \geq T$,
the medians (centers) of the graph are the only nodes with $s_i(t) \geq s_j(t)$ ($h_i(t) \geq h_j(t)$) for all $j \in \Ne{i}$.
\end{theorem}
This theorem implies that there is a time $T$ after which the $s$-values ($h$-values) of the agents do not change.
When the system reaches this time $T$, we say that the $s$-values ($h$-values) have \emph{stabilized}.

While the self-stabilizing graph median and graph center algorithms can be used to identify an optimal leader (an agent whose $s$ or $h$-value is greater than or equal to all of its neighbors), it does not solve the in-network leader selection problem completely. 
We also require that the network has a single leader throughout the execution of the algorithm, not just after the values stabilize.
Our in-network leader selection algorithm leverages the self-stabilizing algorithms above to locate the optimal leader. It also ensures that
the network has a single leader at all times.  The leadership role is transferred from agent to agent.  After the $h$-values or $s$-values stabilize, the leadership role is transferred to an optimal leader in finite time, and this leader remains the leader for all future rounds of the algorithm.

Our self-stabilizing algorithm for in-network leader selection for \LSTV\ is given in Algorithm~\ref{leader.alg}.
One agent is initially selected as leader.  The agents each  have an $s$-value $s_i(t)$ with an arbitrary initial value.
The algorithm executes in synchronous rounds. Each round is divided into two phases.  In the first phase, the agents update their $s_i(t)$ values according to (\ref{medianalg.eq}).  In the second phase, the agent that currently holds the leadership role checks if
it has any neighbors $j$ with $s_j(t) > s_i(t+1)$.  If so, the leader transfers leadership to a neighbor $k$ with maximal $s_k(t)$.

The algorithm for \LSMV\ is nearly identical to Algorithm~\ref{leader.alg}. The only difference is that agents each store $h_i(t)$ instead of $s_i(t)$, and they update this variable according to (\ref{centeralg.eq}).  The leadership transfer phase is identical, with the $h$-values replacing $s$-values. 
Psueodocode for this algorithm is given in Appendix~\ref{maxleadercode.app}.

\subsection{Algorithm Analysis}

We now prove the that Algorithm~\ref{leader.alg} is a self-stabilizing algorithm for in-network leader selection that selects the leader $u$ that 
minimizes $\Terr{\{u\}}$.  
\begin{theorem} \label{leaderT.thm}
Algorithm~\ref{leader.alg} is a self-stabilizing in-network leader selection algorithm for \LSTV.  
\end{theorem}

To prove this theorem, we first introduce some useful definitions and results from~\cite{BGKP99}.
We define a directed graph $G(s) = (V(s), E(s))$ induced by the $s$-values of the agents after these values have stabilized.
The vertex set $V(s)$ is equal to the vertex set $V$ of the original undirected graph $G$.  The edges in $E(s)$ are defined as,
\begin{align*}
&E(s) = \{(i,j)~|~j \in \Ne{i}~\text{and}\\
&~~~~~~~~~~~~~~~(s(j), j)~\text{is lexicographically the largest}\}.
\end{align*}
It has been shown that the undirected underlying graph of $G(s)$ is connected, contains exactly one cycle, and this cycle is of length 2~\cite{BGKP99}.
Let $i$ and $j$ be the nodes that belong to the unique cycle.  Deleting $(i,j)$ and $(j,i)$ from $G(s)$ results in two directed trees,
$T_i(s)$, rooted at $i$ and $T_j(s)$, rooted at $j$.
\begin{theorem}[Proposition 4.3 and Theorem 4.4 from~\cite{BGKP99}] \label{tree.prop}
Each edge in $T_i(s)$ and in $T_j(s)$ is directed from a node to its parent, and if $k$ is a non-leaf node in 
$T_i(s)$ or $T_j(s)$, then $s_k > s_l$ for each child $l$ of $k$. 
\end{theorem}

We now prove Theorem~\ref{leaderT.thm}. \\
\begin{IEEEproof}
It is clear that, since exactly one node is initially the leader, and the leadership role is passed from one node to another,  
there is exactly one leader at any time.  What remains is to show that the algorithm satisfies the convergence and closure properties in Definition~\ref{stable.def}.

The values $s_i(t)$, $i \in V$ are updated according to the algorithm in~(\ref{medianalg.eq}).  As this algorithm is self-stabilizing, 
in  finite time, the $s$-values stabilize; there is a time $T$ after which no $s_i(t)$ changes.  We denote the 
stable value of $s_i(t)$ by $s_i$.
The optimal leaders  for  \LSTV\ are those nodes $i$ such that $s_i \geq s_j$ for all $i \in j$. 

To show that the closure property holds, we must show that, if one of the medians is the leader agent and the $s$-values have stabilized, then this agent remains the leader in all future rounds of the algorithm.  Suppose that in a round $t \geq T$, $i$ is the leader and $i$ is such that $s_i(t) \geq s_j(t)$ for all $j \in \Ne{i}$.  Then, in phase two of the algorithm round, since $i$ does not have any neighbor $j$ with $s_j(t) > s_i(t)$, agent $i$ will not transfer its leadership.  Therefore, the closure property is satisfied.

To show that the convergence property holds, we must show that given any initial values for $s_i(0)$, $i \in V$ and any initial leader assignment, in finite time, the $s$-values will stabilize and the leadership role will be at a median of the graph. Theorem~\ref{selfstable.thm} guarantees that the $s$-values will stabilize in finite time.  

Suppose, after the $s$-values stabilize, a node $u$ has the leadership role, but it is not
a median of the graph.  We show that in finite time, the leadership role will be transferred to a median of the graph.

When a leader transfers leadership to a neighbor, it selects the neighbor with the maximal $s$-value among all its neighbors.
Therefore, leadership is only transferred across edges in the directed graph $G(s)$.   Suppose, without loss of generality,
that the leadership role is at a node $k \neq i$ in the subtree $T_i(s)$.   By the definition of $G(s)$, the neighbor $l$ of $k$ with maximal $s_l$ is the parent of $k$ in the tree.  Further, by Theorem~\ref{tree.prop},  $s_l > s_k$.  Therefore, in the next round, the leadership role will be transferred to agent $l$.  By similar reasoning, in a finite number of rounds, the leadership role will be transferred up the tree until it reaches agent $i$.  
If $s_i > s_j$, then $i$ is the unique median of the graph, since $s_i > s_k$ for all $k \in \Ne{i}$.  
Similarly, if $s_i = s_j$, then both $i$ and $j$ are medians of the graph.  In either of these cases, the leadership role has reached a median.  
If $s_j > s_i$, $j$ is the unique median of the graph, since $s_j > s_k$ for all $k \in \Ne{i}$.
Since $s_i < s_j$ and all other neighbors $k$ of $i$ are such that $s_k < s_i$, in the next round, agent $i$ will transfer leadership to agent $j$.
\end{IEEEproof}

A similar result also holds for the in-network leader selection algorithm that  selects the leader $u$ that 
minimizes $\Merr{\{u\}}$.  We omit this proof since it is similar to the proof of Theorem~\ref{leaderT.thm}.
\begin{theorem}
Algorithm~\ref{leaderM.alg} is a self-stabilizing in-network leader selection algorithm for \LSMV.  
\end{theorem}
  
For  Algorithm~\ref{leader.alg}, 
the number of rounds until the $s$-values stabilize, starting from any initial state, is in $\Theta(d)$, where $d$ is the maximum distance from 
the edge of the network (the nodes $i$ with $| \Ne{i}| = 1$) to a median.
Once the $s$-values stabilize, it takes at most $d$ rounds for the leadership role to be transferred to the median node. Therefore
the running time of Algorithm~\ref{leader.alg} is in $\Theta(d)$.

For Algorithm~\ref{leaderM.alg}, the number of rounds until the $h$-values stabilize is $\Theta(r)$ where $r$ is the radius of the graph. Once the $h$-values stabilize, it takes at most $r$  rounds for the leadership role to be transferred to the center.  Therefore the running time of Algorithm~\ref{leaderM.alg} is in $\Theta(r)$.

%%%%%%%%%%%%%%%%%%%%%%%%%%%%%%%%%%%%%%%%%%%%%%%%%%%%%%%%%%%%%%%%%%%%%%%%%
%%%%%%%%%%%%%%%%%%%%%%%%%%%%%%%%%%%%%%%%%%%%%%%%%%%%%%%%%%%%%%%%%%%%%%%%%
%%%%%%%%%%%%%%%%%%%%%%%%%%%%%%%%%%%%%%%%%%%%%%%%%%%%%%%%%%%%%%%%%%%%%%%%%
%%%%%%%%%%%%%%%%%%%%%%%%%%%%%%%%%%%%%%%%%%%%%%%%%%%%%%%%%%%%%%%%%%%%%%%%%
%%%%%%%%%%%%%%%%%%%%%%%%%%%%%%%%%%%%%%%%%%%%%%%%%%%%%%%%%%%%%%%%%%%%%%%%%

\vspace{-.2cm}

\section{Extension to Weighted Graphs}
\label{extend.sec}

Several recent works have investigated a generalization of the $k$-leader selection problem where the objective is for every agent
 to maintain specified differences between its state and the states of its neighbors,
\begin{equation} \label{diff.eq}
x_i(t) - x_j(t) = \pos{i}{j}~~~~\text{for all}~(i,j) \in E,
\end{equation}
 where $\pos{i}{j}$ denotes the desired difference, for example, the position of node $i$ relative to node $j$~\cite{BH06,BH08,CBP14}.

The states of the leaders are the reference states. Let $\hat{x}$ denote vector of desired states that satisfy (\ref{diff.eq}) when the leader states are fixed.
 A follower  updates its state based on noisy measurements of the differences between its state and the states
of its neighbors. The dynamics of each follower agent is given by
\[
\dot{x}_i(t) = -\sum_{j \in \Ne{i}} W_{ij}  \left(x_i(t) - x_j(t) - \pos{i}{j}+ \epsilon_{ij}(t)\right),
\]
where $W_{ij}$ is the weight for link $(i,j)$ and
 $\epsilon_{ij}(t)$, $(i,j) \in E$, are independent, identically distributed zero-mean white noise processes with autocorrelation functions 
$\expec{\epsilon_{ij}(t)\epsilon_{ij}(t + \tau)} = \nu_{ij} \delta(\tau)$.
The link weights are chosen so that $W_{ij} = \frac{\nu_{ij}}{\Deg{i}}$, where $\Deg{i} = \sum_{j \in \Ne{i}}\frac{1}{\nu_{ij}}$, which correspond to the best linear unbiased estimator of the leader agents state when $x_j(t) = \hat{x}_{j}$ for all $j \in \Ne{k}$~\cite{BH06}.

In this system, the leader-follower dynamics are depend on a weighted Laplacian matrix, also called the conductance matrix,
which defined as
\[
L_{ij} = \left\{ \begin{array}{ll} 
- \frac{1}{\nu_{ij}} & ~\text{if}~(i,j) \in E \\
\Deg{i} &~\text{if}~i=j \\
0 &~\text{otherwise.}
\end{array} \right.
\]
Note that if  $\nu_{ij} = 1$ for all $(i,j) \in E$, then $L$ is standard, unweighted Laplacian matrix. 
We assume that the weights are symmetric, meaning $\nu_{ij} = \nu_{ji}$ for all $(i,j) \in E$.

Let $\sigma_i$ be the steady-state variance of the deviation from $\hat{x}$,
\[
\sigma_i = \lim_{t \rightarrow \infty} \expec{x_i(t) - \hat{x}_i}^2.
\]
As with the unweighted Laplacian, $\sigma_i = \frac{1}{2} (\Lff^{-1})_{ii}$, where $\Lff$ is the follower-follower submatrix
of the weighted Laplacian.
It has also been shown that $\sigma_i = \frac{1}{2}\resistwt{i}{j}$, where $\resistwt{i}{j}$ is the resistance distance between 
$i$ and $j$ in the electrical network
where each edge $(i,j)$ has a $\nu_{ij}$ resistor~\cite{BH06}.  

If we define $\fdist{i}{j} = \resistwt{i}{j}$, then, as before, a solution to the $p$-median problem for $p=1$ is a solution to \LSTV,
and a solution to the $p$-center problem for $p=1$, is a solution to \LSMV. It is also straightforward to extend our in-network 
leader selection algorithms to this weighted graph setting.

Let $\Gwt=(V,\Ewt)$ be a weighted, connected acyclic graph, where the weight of edge $(i,j)$ is $c_{ij} = 1/\nu_{ij}$.
In this graph, the resistance distance between nodes $i$ and $j$ is equal to the graph distance, where the graph distance
is the sum of the weights of the edges in the path between nodes $i$ and $j$~\cite{KR93}.  
Our in-network leader selection algorithms must find the median and center for this weighted graph.

We first state an important theorem about the medians of acyclic graphs with positive edge weights
\begin{theorem}[Lemma 7.2 in \cite{BGKP99}] The medians of an acyclic, connected graph remain unchanged independent of any change in the weight of the edges.
\end{theorem} 
Since the location of the graph medians, and hence the identity of an optimal leader, do not depend on the edge weights,
Algorithm~\ref{leader.alg} also solves the in-network leader selection for \LSTV\ in acyclic weighted graphs.

It has been shown that the self-stabilizing center finding algorithm in (\ref{centeralg.eq}) can be extended to weighed graphs by making small modifications~\cite{BGKP99}. Let $\Nhminusw{i}$ be defined as follows,
\begin{align*}
\Nhminusw{i} &= \{h_j(t) + c_{ij}~|~(i,j) \in \Ne{i}\} \\
&~~~~~~~~~~~~~~- \max\{h_j(t) + c_{ij}~|~(i,j) \in \Ne{i}\}.
\end{align*}
The modified center-finding algorithm is as follows,
\[
h_i(t+1) = \left\{ \begin{array}{ll}
0 &~\text{if}~|\Ne{i}| = 1 \\
1 + \max \left(\Nhminusw{i}\right)&~\text{otherwise}.
\end{array}\right.
\]

By incorporating these small changes to the update to the $h$-values, Algorithm~\ref{leaderM.alg} can be used to solve the in-network leader selection problem for 
\LSMV\ in acyclic weighted graphs.

\section{Conclusion}
\label{conclusion.sec}
In this work, we have posed a new leader selection problem called the in-network leader selection problem,
whereby agents must collaborate to find the leader that optimizes a chosen performance measure.
We have considered two performance measures, the total steady-state variance of the deviation and the maximal 
steady-state variance of the deviation.
We have shown that finding a leader that minimizes the total variance is equivalent to solving the $p$-median facility location problem for $p=1$
and that finding a leader that minimizes the maximal variance is equivalent to the $p$-center facility location problem for $p=1$.
Leveraging the connections to these problems, we have developed two self-stabilizing in-network leader selection algorithms, one for 
each performance measure.
Finally, we have shown how our algorithms can be extended to weighted graphs.
In future work, we plan to investigate generalizing our approach to in-network algorithms for the $k$-leader selection problem where $k$ is greater than one.


%\section*{Acknowledgment}

\begin{algorithm}
\caption{Algorithm for in-network leader selection for \LSMV.} \label{leaderM.alg}
\begin{algorithmic}
\State {Send $h_i(t)$ to all  agents $j$ in $\Ne{i}$}
\State {Receive $h_j(t)$ from all $j \in \Ne{i}$}
\If {$|\Ne{i}| = 1$}
	\State{$h_i(t+1) \gets 0$}
\Else
	\State {$h_i(t+1)  \gets  1+ \max \left(\Nhminus{i}\right)$}
\EndIf
\If {(\textit{leader} = TRUE) and \\~~~~~~~~~~~($h_j(t) > h_i(t+1)$ for some $j \in \Ne{i}$) }
	\State{$k =  \arg \max_{j \in \Ne{i}} h_j(t)$}
	\State{$\textit{leader} \gets$  FALSE}
	\State{transfer leadership to agent $k$}
\ElsIf{receive leadership transfer from neighbor}
	\State{$\textit{leader} \gets$ TRUE}
\EndIf
\end{algorithmic}
\end{algorithm}

\appendices

\section{Illustration that $\Merr{\cdot}$ is not Super-Modular} \label{notsubmod.app}
\begin{figure}
\centering
\includegraphics[scale=.5]{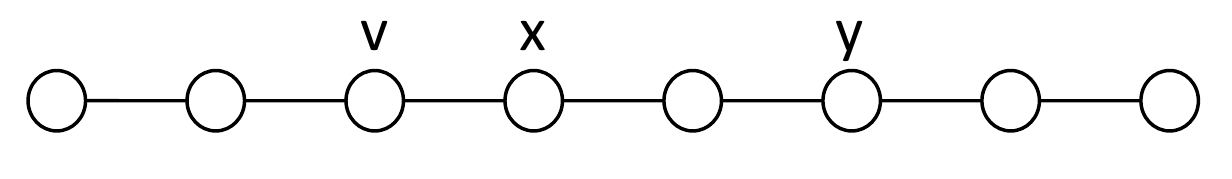}
\caption{Example network demonstrating that the error measure $\Merr{\cdot}$ is not super-modular.
Let $A =\{x\}$ and $B=\{x,y\}$.   Then,  $\Merr{A} = 2$,  $\Merr{A \cup \{v\}} =  2$, $\Merr{B} = 1.5$, and $\Merr{B \cup \{v\}} = 1$. 
Thus, $\Merr{A} - \Merr{A \cup \{v\}} < \Merr{B} - \Merr{B \cup \{v\}}$.}
 \label{line.fig}
\end{figure}

A super-modular function is defined as follows.  Let $V$ be a finite set and let $A$ and $B$ be sets with $A \subseteq B \subseteq V$.
A function $f$ is \emph{super-modular} if and only if  for all  $v \in V \setminus B$,
\[
f(A) - f(A \cup \{v\}) \geq f(B) -  f(B \cup \{v\}).
\] 
In Figure~\ref{line.fig}, we give an example that illustrates that $\Merr{\cdot}$ is not a super-modular function.  
\newpage
\section{Pseudocode for In-Network Leader Selection Algorithm for \LSMV} \label{maxleadercode.app}

Pseudocode for the self-stabilizing algorithm that solves \LSMV\ is given in Algorithm~\ref{leaderM.alg}.
\bibliographystyle{IEEEtran}
\bibliography{leader}

\begin{thebibliography}{10}
\providecommand{\url}[1]{#1}
\csname url@samestyle\endcsname
\providecommand{\newblock}{\relax}
\providecommand{\bibinfo}[2]{#2}
\providecommand{\BIBentrySTDinterwordspacing}{\spaceskip=0pt\relax}
\providecommand{\BIBentryALTinterwordstretchfactor}{4}
\providecommand{\BIBentryALTinterwordspacing}{\spaceskip=\fontdimen2\font plus
\BIBentryALTinterwordstretchfactor\fontdimen3\font minus
  \fontdimen4\font\relax}
\providecommand{\BIBforeignlanguage}[2]{{%
\expandafter\ifx\csname l@#1\endcsname\relax
\typeout{** WARNING: IEEEtran.bst: No hyphenation pattern has been}%
\typeout{** loaded for the language `#1'. Using the pattern for}%
\typeout{** the default language instead.}%
\else
\language=\csname l@#1\endcsname
\fi
#2}}
\providecommand{\BIBdecl}{\relax}
\BIBdecl

\bibitem{RBM05}
W.~Ren, R.~W. Beard, and T.~W. McLain, ``Coordination variables and consensus
  building in multiple vehicle systems,'' \emph{Cooperative Control}, vol. 309,
  pp. 171--188, 2005.

\bibitem{EKPS04}
J.~Elson, R.~M. Karp, C.~H. Papadimitriou, and S.~Shenker, ``Global
  synchronization in sensornets,'' in \emph{Latin American Theoretical
  Informatics}, 2004, pp. 609--624.

\bibitem{BH09}
P.~Barooah and J.~Hespanha, ``Error scaling laws for linear optimal estimation
  from relative measurements,'' \emph{IEEE Transactions on Information Theory},
  vol.~55, no.~12, pp. 5661--5673, Dec 2009.

\bibitem{BH06}
------, ``Graph effective resistance and distributed control: Spectral
  properties and applications,'' in \emph{Proc. 45th IEEE Conf. on Decision and
  Control}, Dec 2006, pp. 3479=--3485.

\bibitem{PB10}
S.~Patterson and B.~Bamieh, ``Leader selection for optimal network coherence,''
  in \emph{Proc. 49th IEEE Conf. on Decision and Control}, 2010, pp.
  2692--2697.

\bibitem{LFJ14}
F.~Lin, M.~Fardad, and M.~Jovanovic, ``Algorithms for leader selection in
  stochastically forced consensus networks,'' \emph{IEEE Transactions on
  Automatic Control}, vol.~59, no.~7, pp. 1789--1802, Jul 2014.

\bibitem{CBP14}
A.~Clark, L.~Bushnell, and R.~Poovendran, ``A supermodular optimization
  framework for leader selection under link noise in linear multi-agent
  systems,'' \emph{IEEE Transactions on Automatic Control}, vol.~59, no.~2, pp.
  283--296, Feb 2014.

\bibitem{N78}
G.~L. Nemhauser, L.~A. Wolsey, and M.~L. Fisher, ``An analysis of
  approximations for maximizing submodular set functions,'' \emph{Mathematical
  Programming}, vol.~14, no.~1, pp. 265--294, 1978.
\newpage
\bibitem{CABP14}
A.~Clark, B.~Alomair, L.~Bushnell, and R.~Poovendran, ``Minimizing convergence
  error in multi-agent systems via leader selection: A supermodular
  optimization approach,'' \emph{IEEE Transactions on Automatic Control},
  vol.~59, no.~6, pp. 1480--1494, Jun 2014.


\bibitem{FL13}
K.~Fitch and N.~Leonard, ``Information centrality and optimal leader selection
  in noisy networks,'' in \emph{Proc. 52nd IEEE Conf. on Decision and Control},
  Dec 2013, pp. 7510--7515.

\bibitem{HD02}
H.~W. Hamacher and Z.~Drezner, \emph{Facility location: applications and
  theory}.\hskip 1em plus 0.5em minus 0.4em\relax Springer, 2002.

\bibitem{GBS08}
A.~Ghosh, S.~Boyd, and A.~Saberi, ``Minimizing effective resistance of a
  graph,'' \emph{SIAM Review}, vol.~50, no.~1, pp. 37--66, Feb 2008.

\bibitem{KR93}
D.~Klein and M.~Randic, ``Resistance distance,'' \emph{Journal of Mathematical
  Chemistry}, vol.~12, no.~1, pp. 81--95, 1993.

\bibitem{BJMP12}
B.~Bamieh, M.~Jovanovic, P.~Mitra, and S.~Patterson, ``Coherence in large-scale
  networks: Dimension-dependent limitations of local feedback,'' \emph{IEEE
  Transactions on Automatic Control}, vol.~57, no.~9, pp. 2235--2249, Sep 2012.

\bibitem{D74}
E.~W. Dijkstra, ``Self-stabilizing systems in spite of distributed control,''
  \emph{Communications of the ACM}, vol.~17, no.~11, pp. 643--644, Nov 1974.

\bibitem{L96}
N.~Lynch, \emph{Distributed Algorithms}.\hskip 1em plus 0.5em minus 0.4em\relax
  USA: Morgan Kaufmann Publishers, Inc., 1996.

\bibitem{KH79}
O.~Kariv and S.~L. Hakimi, ``An algorithmic approach to network location
  problems. ii: The p-medians,'' \emph{SIAM Journal on Applied Mathematics},
  vol.~37, no.~3, pp. 539--560, 1979.

\bibitem{BGKP99}
S.~C. Bruell, S.~Ghosh, M.~H. Karaata, and S.~V. Pemmaraju, ``Self-stabilizing
  algorithms for finding centers and medians of trees,'' \emph{SIAM Journal on
  Computing}, vol.~29, no.~2, pp. 600--614, 1999.

\bibitem{BH08}
P.~Barooah and J.~Hespanha, ``Estimation from relative measurements: Electrical
  analogy and large graphs,'' \emph{IEEE Transactions on Signal Processing},
  vol.~56, no.~6, pp. 2181--2193, Jun 2008.

\end{thebibliography}

\end{document}